\documentclass[12pt,reqno]{amsart}

\usepackage{fullpage,amsmath,amssymb,amsfonts,mathrsfs,bm,graphicx,hyperref}
\usepackage{tikz}

\def\bbone{{\mathchoice {\rm 1\mskip-4mu l} {\rm 1\mskip-4mu l}
{\rm 1\mskip-4.5mu l} {\rm 1\mskip-5mu l}}}

\newtheorem{theorem}{Theorem}[section]
\newtheorem{remark}{Remark}[section]

\newtheorem{lemma}{Lemma}[section]
\newtheorem{proposition}{Proposition}[section]

\begin{document}

\author{Abdelmalek Abdesselam}
\address{Abdelmalek Abdesselam, Department of Mathematics,
P. O. Box 400137,
University of Virginia,
Charlottesville, VA 22904-4137, USA}
\email{malek@virginia.edu}

\author{Shannon Starr}
\address{Shannon Starr, Department of Mathematics,
University of Alabama at Birmingham,
1402 10th Avenue South, Birmingham, AL 35294-1241, USA}
\email{slstarr@uab.edu}

\title{A central limit theorem for a generalization of the Ewens measure to random tuples of commuting permutations}

\begin{abstract}
We prove a central limit theorem (CLT) for the number of joint orbits of random tuples of commuting permutations. In the uniform sampling case this generalizes the classic CLT of Goncharov for the number of cycles of a single random permutation. We also consider the case where tuples are weighted by a factor other than one, per joint orbit. We view this as an analogue of the Ewens measure, for tuples of commuting permutations, where our CLT generalizes the CLT by Hansen. Our proof uses saddle point analysis, in a context related to the Hardy-Ramanujan asymptotics and the theorem of Meinardus, but concerns a multiple pole situation.
The proof is written in a self-contained manner, and hopefully in a manner accessible to a wider audience. We also indicate several open directions of further study related to probability, combinatorics, number theory, an elusive theory of random commuting matrices, and perhaps also geometric group theory.
\end{abstract}

\maketitle


\section{Introduction}

\subsection{The main result}

For $n\ge 0$, the symmetric group $\mathfrak{S}_n$ comes with a natural left action $(\sigma,i)\mapsto \sigma(i)$ on the set $[n]:=\{1,2,\ldots,n\}$. If $H$ is a subgroup of $\mathfrak{S}_n$, then the restriction of this action to $H$ gives us access to the enumerative quantity $\kappa_n(H):=|[n]/H|$, where we used $|\cdot|$ to denote the cardinality of finite sets, and where $[n]/H$ is the set of orbits for the action of $H$ on $[n]$. We will be particularly interested in the situation where $H=\langle\sigma_1,\ldots,\sigma_{\ell}\rangle$ is the subgroup generated by a random ordered $\ell$-tuple $(\sigma_1,\ldots,\sigma_{\ell})$ of pairwise commuting permutations, for some fixed $\ell\ge 1$. Our main result is a central limit theorem for the number of joint orbits $\kappa_n(\langle\sigma_1,\ldots,\sigma_{\ell}\rangle)$ or simply $\kappa_{\ell,n}(\sigma_1,\ldots,\sigma_{\ell})$. 
Based on a description in \cite{AbdesselamBDV}, the tuples with a single orbit are twisted discrete $\ell$-dimensional tori,
and these form the indivisible units for tuples with $k>1$ orbits.
See also Section \ref{subsec:relations}.

For $\ell=2$ if $n=d_1d_2$ then there are pairs of permutations $(\sigma_1,\sigma_2)$ corresponding
to a discrete $d_1\times d_2$ torus with a single twist $t \in \{0,\dots,d_2-1\}$.
Examples for $d_1=d_2=4$ and twist $t \in \{0,1,2,3\}$ are shown in Figure \ref{fig:torus}.
\begin{figure}
\begin{center}
\begin{tikzpicture}[xscale=0.65,yscale=0.65]
\begin{scope}[xshift=-6cm]
\foreach \x in {0,1,2,3}
	{\draw (\x,5) -- +(0,1);}
\foreach \y in {0,1,2,3,4,5}
	{\draw (0,\y) -- (4,\y);
	 \fill[white] (4,\y) circle (2mm);
	 \draw (4,\y) circle (1mm);}
\foreach \x in {0,1,2,3}
	{	\draw (\x,0) -- (\x,5);
		\foreach \y in {0,1,2,3,4,5}
			{   \fill[white] (\x,\y) circle (2mm);
				\fill (\x,\y) circle (1mm);
			}
		\fill[white] (\x,6) circle (2mm);
		\draw (\x,6) circle (1mm);
	}
\end{scope}
\draw (3,5) -- (0,6);
\foreach \x in {0,1,2}
	{\draw[line width=5pt,white] (\x,5) -- +(1,1);
	\draw (\x,5) -- +(1,1);}
\foreach \y in {0,1,2,3,4,5}
	{\draw (0,\y) -- (4,\y);
	 \fill[white] (4,\y) circle (2mm);
	 \draw (4,\y) circle (1mm);}
\foreach \x in {0,1,2,3}
	{	\draw (\x,0) -- (\x,5);
		\foreach \y in {0,1,2,3,4,5}
			{   \fill[white] (\x,\y) circle (2mm);
				\fill (\x,\y) circle (1mm);
			}
		\fill[white] (\x,6) circle (2mm);
		\draw (\x,6) circle (1mm);
	}
\begin{scope}[xshift=6cm]
\draw (2,5) -- (0,6);
\draw (3,5) -- (1,6);
\foreach \x in {0,1}
	{\draw[line width=5pt,white] (\x,5) -- +(2,1);
	  \draw (\x,5) -- +(2,1);}
\foreach \y in {0,1,2,3,4,5}
	{\draw (0,\y) -- (4,\y);
	 \fill[white] (4,\y) circle (2mm);
	 \draw (4,\y) circle (1mm);}
\foreach \x in {0,1,2,3}
	{	\draw (\x,0) -- (\x,5);
		\foreach \y in {0,1,2,3,4,5}
			{   \fill[white] (\x,\y) circle (2mm);
				\fill (\x,\y) circle (1mm);
			}
		\fill[white] (\x,6) circle (2mm);
		\draw (\x,6) circle (1mm);
	}
\end{scope}
\begin{scope}[xshift=12cm]
\draw (1,5) -- (0,6);
\draw (2,5) -- (1,6);
\draw (3,5) -- (2,6);
\draw[line width=5pt,white] (0,5) -- +(3,1);
\draw (0,5) -- +(3,1);
\foreach \y in {0,1,2,3,4,5}
	{\draw (0,\y) -- (4,\y);
	 \fill[white] (4,\y) circle (2mm);
	 \draw (4,\y) circle (1mm);}
\foreach \x in {0,1,2,3}
	{	\draw (\x,0) -- (\x,5);
		\foreach \y in {0,1,2,3,4,5}
			{   \fill[white] (\x,\y) circle (2mm);
				\fill (\x,\y) circle (1mm);
			}
		\fill[white] (\x,6) circle (2mm);
		\draw (\x,6) circle (1mm);
	}
\end{scope}
\end{tikzpicture}
\end{center}
\caption{Four discrete $\ell=2$ dimensional tori with a twist. 
The white circles represent the vertices on the opposite face, repeated so as to 
draw the graphs without extraneous crossings.
\label{fig:torus}}
\end{figure}
It helps to consider permutations of the numbers $\{0,\dots,n-1\}$ instead of $[n]$ for describing these discrete tori with twists.
We give canonical examples, and all other examples are obtained by conjugation.
The first permutation can be written as a product of cycles
$\sigma_1 = \rho_0\rho_1\cdots \rho_{d_1-1}$ where $\rho_k = (kd_2,kd_2+1\dots,kd_2+d_2-1)$ for $k\in\{0,\dots,d_1-1\}$.
Then $\sigma_2$ may be viewed as the bijection $\varphi$ from $\{0,\dots,n-1\}$ to itself, as follows.
Suppose $i=j+kd_2$ for $k \in \{0,\dots,d_1-1\}$ and $j \in \{0,\dots,d_2-1\}$.
Then $\varphi(i)$ has two cases: if $k+1<d_1$ the $\varphi(i) = j+(k+1)d_2$;
but if $k=d_1-1$ then $\varphi(i) = \operatorname{mod}(j+t,d_2)$ where $t \in \{0,\dots,d_2-1\}$ determines the single twist 
for $\ell=2$.

To state our results more precisely, let us introduce the set of commuting tuples of permutations
\[
\mathcal{C}_{\ell,n}:=\{(\sigma_1,\ldots,\sigma_{\ell})\in \mathfrak{S}_n^{\ell}\ |\ \forall i,j,\ \sigma_i\sigma_j=\sigma_j\sigma_i\}\ ,
\]
and the numbers
\[
A(\ell,n,k):=\left|\left\{
(\sigma_1,\ldots,\sigma_{\ell})\in \mathcal{C}_{\ell,n}\ |\ 
\kappa_{\ell,n}(\sigma_1,\ldots,\sigma_{\ell})=k
\right\}\right|\ ,
\]
for $0\le k\le n$, as well as the polynomials
\[
\mathcal{H}_{\ell,n}(x):=\frac{1}{n!}\sum_{k=0}^{n}A(\ell,n,k)\ x^k\ .
\] 
When $x>0$, we define the probability measure $\mathbb{P}_{\ell,n,x}$ on $\mathcal{C}_{\ell,n}$ given by weighting a tuple $(\sigma_1,\ldots,\sigma_{\ell})$ by $x^{\kappa_{\ell,n}(\sigma_1,\ldots,\sigma_{\ell})}$.
Namely, the corresponding probability mass function (PMF) is
\[
\mathbb{P}_{\ell,n,x}(\left\{(\sigma_1,\ldots,\sigma_{\ell})\right\})=
\frac{x^{\kappa_{\ell,n}(\sigma_1,\ldots,\sigma_{\ell})}}{n!\ \mathcal{H}_{\ell,n}(x)}\ .
\]
We view $\mathbb{P}_{\ell,n,x}$ as a natural generalization of the Ewens measure with parameter $x$ to tuples of commuting permutations. Indeed, when $\ell=1$, 
$\kappa_{1,n}(\sigma)$ is just the number of cycles of the single permutation $\sigma$, and
the well known expansion of the Pochhammer symbol in terms of Stirling numbers of the first kind gives the following specialization of the above PMF formula
\[
\mathbb{P}_{1,n,x}(\left\{\sigma\right\})=
\frac{x^{\kappa_{1,n}(\sigma)}}{x(x+1)\cdots(x+n-1)}\ .
\]
This is the PMF of the much studied Ewens distribution on the symmetric group, which was introduced in the applied context of mathematical biology~\cite{Ewens}, but is also
an important notion from a pure mathematics point of view (see~\cite{Olshanski,Crane,Tavare} for insightful surveys, and the book~\cite{ArratiaBTbook} for an in-depth introduction to the Ewens measure and many other related topics in probability theory and combinatorics). The uniform $x=1$ case was of course studied much earlier than the article by Ewens, but we also note that the particular weight $x=2$ per cycle (in a more difficult spatially dependent rather than mean field situation) features in the work of T\'oth on quantum spin systems~\cite[\S5]{Toth}.

When the random $\ell$-tuple $(\sigma_1,\ldots,\sigma_{\ell})\in\mathcal{C}_{\ell,n}$ is sampled according to the measure $\mathbb{P}_{\ell,n,x}$,
this gives rise to the random variable $\mathsf{K}_{\ell,n}:=\kappa_{\ell,n}(\sigma_1,\ldots,\sigma_{\ell})$. Before stating our CLT for $\mathsf{K}_{\ell,n}$, as $n\rightarrow \infty$, let us define for notational convenience the constant
\begin{equation}
\mathcal{K}_{\ell}:=(\ell-1)!\ \zeta(2)\zeta(3)\cdots\zeta(\ell)\ ,
\label{Kdefeq}
\end{equation}
which features special values of the Riemann zeta function, and which reduces to $\mathcal{K}_{1}=1$ when $\ell=1$ since the products are empty. Note that when we write an asymptotic equivalence $f(u)\sim g(u)$ when the argument $u$ goes to some limit, we mean the precise statement $\lim \frac{f(u)}{g(u)}=1$, as usual in asymptotic analysis. Our use of the Landau symbols $o(\cdot)$ and $O(\cdot)$ also follows the standard custom (see, e.g.,~\cite{deBruijn}).

\begin{theorem}\label{mainthm}
For any $\ell\ge 2$, and any $x>0$, as $n\rightarrow\infty$, 
the leading asymptotics of the mean and variance of the random variables $\mathsf{K}_{\ell,n}$ are given by
\begin{eqnarray}
\mathbb{E}\mathsf{K}_{\ell,n}&\sim & 
\frac{(x\mathcal{K}_{\ell})^{\frac{1}{\ell}}}{\ell-1}\times n^{\frac{\ell-1}{\ell}}\ , 
\label{Easymeq} \\
{\rm Var}(\mathsf{K}_{\ell,n}) & \sim & 
\frac{(x\mathcal{K}_{\ell})^{\frac{1}{\ell}}}{\ell(\ell-1)}\times n^{\frac{\ell-1}{\ell}}
\ . \label{Vasymeq}
\end{eqnarray}
Moreover, the normalized random variables
\[
\frac{\mathsf{K}_{\ell,n}-\mathbb{E}\mathsf{K}_{\ell,n}}{\sqrt{{\rm Var}(\mathsf{K}_{\ell,n})}}
\]
converge in distribution and in the sense of moments to the standard Gaussian $\mathcal{N}(0,1)$. Namely, we have
\[
\lim\limits_{n\rightarrow\infty}\mathbb{E}\left[f\left(
\frac{\mathsf{K}_{\ell,n}-\mathbb{E}\mathsf{K}_{\ell,n}}{\sqrt{{\rm Var}(\mathsf{K}_{\ell,n})}}
\right)\right]=\frac{1}{\sqrt{2\pi}}\int\limits_{-\infty}^{\infty}f(s)\ e^{-\frac{s^2}{2}}\ {\rm d}s\ ,
\]
for all $f$'s that are bounded continuous functions, or polynomials.
\end{theorem}

\begin{figure}
\begin{center}
\begin{tikzpicture}
\draw  (0,0) node[] {\includegraphics[width=13cm]{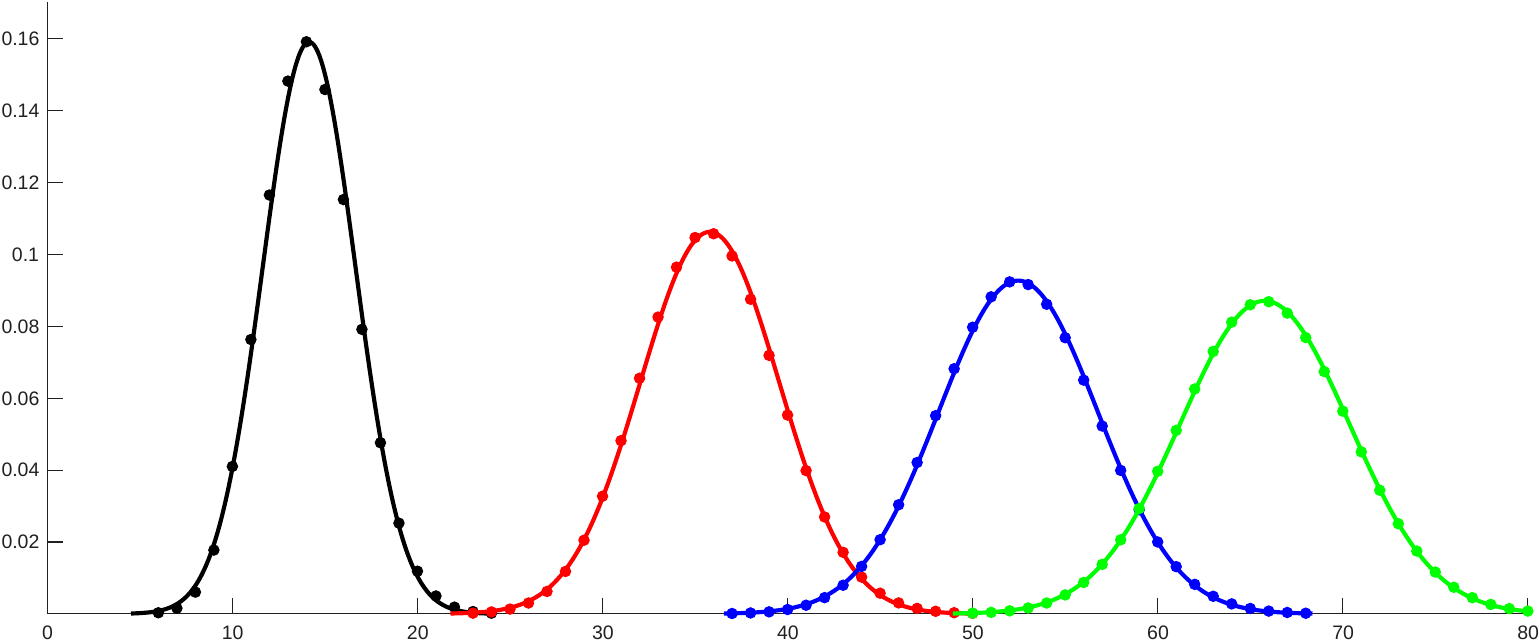}};
\end{tikzpicture}
\end{center}
\caption{Plots of the probability mass functions for $\mathsf{K}_{\ell,n}$ for $\ell=2$, $n=150$ and choices of $x$:
$x=1$ in black, $x=9$ in red, $x=25$ in blue and $x=49$ in green.
The pmf is the dotted plot. The solid curves are the Gaussian probability density functions with the same mean and variance, so that one may compare
the agreement. (We have applied a cutoff for very small values of probability by starting the window at $10^{-4}$ on the vertical
scale.)
\label{fig:CLT}}
\end{figure}
In Figure \ref{fig:CLT} we have plotted the probability mass function for several values of $x$, along with the Gaussian approximation.
The data for this image is given in Appendix \ref{sec:data}.

Our statement does not cover the $\ell=1$ case, but the CLT is also true in the latter, provided one interprets $\frac{1}{\ell-1}n^{\frac{\ell-1}{\ell}}$ as $\ln n$. 
Indeed, one has the convergence in distribution
\[
\frac{\mathsf{K}_{1,n}-\mathbb{E}\mathsf{K}_{1,n}}{\sqrt{{\rm Var}(\mathsf{K}_{1,n})}}\Longrightarrow \mathcal{N}(0,1)\ ,
\]
proved by Hansen as the time 1 projection of a functional CLT for a process related to the Ewens measure~\cite{Hansen}. For $\ell=1$, the mean and variance are both asymptotically equivalent to $x\ln n$.
The $x=1$ uniform sampling case is much older and due to Goncharov~\cite{Goncharov1,Goncharov2} and is often given as an instructive example in graduate probability textbooks (see, e.g.,~\cite[Example 27.3]{Billingsley}).
Using the Feller coupling idea~\cite[p. 815]{Feller}, one can express $\mathsf{K}_{1,n}$ as a sum of $n$ independent Bernoulli random variables of parameters $\frac{1}{n-j+1}$, $1\le j\le n$. Since these are not identically distributed, they provide a nice pedagogical example of application of the Lindeberg CLT.
We also note the remarkable continuity of the formulas for the asymptotic mean and variance of the $\mathsf{K}_{\ell,n}$ as one varies $\ell$ over the full range $\ell\ge 1$, with the above caveat of logarithmic interpretation of the power law when $\ell=1$.
This is reminiscent of the behavior of models of statistical mechanics at the upper critical dimension (see~\cite{BauerschmidtBS} for a general introduction).

\subsection{A brief outline of the proof}

The only tool from probability theory we will use is the following continuity theorem regarding convergence of moment generating functions on the real line.
\begin{theorem}\label{Curtissthm}
Let $\mathsf{X}_n$ be a sequence of real random variables, and suppose $\exists s_0>0$, such that $\forall s\in(-s_0,s_0)$, the moment generating functions
\[
M_n(s):=\mathbb{E}\left[e^{s\mathsf{X}_n}\right]
\]
are well defined (the integrals converge), and have pointwise limits
\[
M(s):=\lim\limits_{n\rightarrow \infty}M_n(s)\ .
\]
Then there exists a unique random variable $\mathsf{X}$ (or distribution rather), such that $\mathsf{X}_n$ converge in distribution to $\mathsf{X}$. The moment generating function of $\mathsf{X}$ is well defined on the interval $(-s_0,s_0)$ and coincides with $M(s)$. Moreover, the convergence also holds in the sense of moments, i.e., $\lim_{n\rightarrow\infty}\mathbb{E}[\mathsf{X}_{n}^{p}]=\mathbb{E}[\mathsf{X}^p]$ for every nonnegative integer $p$. 
\end{theorem}
This is a classical theorem of Curtiss~\cite{Curtiss}. Although the part in the conclusion about the convergence of moments is not often mentioned explicitly, it easily follows by a standard integration to the limit argument (see Theorem 25.12 and its corollary in~\cite{Billingsley}).
The key proposition needed for the proof of Theorem \ref{mainthm} is the following.
\begin{proposition}\label{keyprop}
Suppose $\ell\ge 2$ and $x>0$, and consider the previous random variables $\mathsf{K}_{\ell,n}$ with distribution determined by the measure $\mathbb{P}_{\ell,n,x}$. 
Let $(a_n)$ be a sequence of real numbers, and let $(b_n)$ be a sequence of positive real numbers such that, as $n\rightarrow\infty$,
\begin{eqnarray*}
a_n &=& 
x Z_{\ell-1}^{[\ell]}\left(
\left(Z_{\ell}^{[\ell]}\right)^{-1}\left(\frac{n}{x}\right)
\right)
+o\left(n^{\frac{\ell-1}{2\ell}}\right)\ ,  \\
b_n & = & 
\frac{(x\mathcal{K}_{\ell})^{\frac{1}{2\ell}}}{\sqrt{\ell(\ell-1)}}\times n^{\frac{\ell-1}{2\ell}}\times (1+o(1))\ ,
\end{eqnarray*}
where $Z_{\ell-1}^{[\ell]}$ and $Z_{\ell}^{[\ell]}$ are explicit bijective functions $(0,\infty)\rightarrow(0,\infty)$ defined further below.
Then, for all $s\in\mathbb{R}$, we have
\[
\lim\limits_{n\rightarrow\infty}\ln
\mathbb{E}\left[\exp\left(s\left(\frac{\mathsf{K}_{\ell,n}-a_n}{b_n}\right)\right)\right]
=\frac{s^2}{2}\ .
\]
\end{proposition}
By doing several passes applying Proposition \ref{keyprop} and Theorem \ref{Curtissthm} to suitable sequences $a_n$ and $b_n$, we deduce our CLT in the clean form given in Theorem \ref{mainthm}. As for establishing Proposition \ref{keyprop}, this is done using a saddle point analysis in the spirit of the book~\cite{FlajoletS}.
For any fixed $x>0$, the ordinary generating function of the sequence $(\mathcal{H}_{\ell,n}(x))_{n\ge 0}$ is given by
\begin{equation}
\mathcal{G}_{\ell}(x,z):=\sum_{n=0}^{\infty}\mathcal{H}_{\ell,n}(x)\ z^n
=\prod_{\delta_1,\ldots,\delta_{\ell-1}=1}^{\infty}
\left(1-z^{\delta_1\cdots\delta_{\ell-1}}\right)^{-x\delta_1^{\ell-2}\delta_2^{\ell-3}\cdots\delta_{\ell-2}}\ ,
\label{BryanFeq}
\end{equation}
thanks to a formula by Bryan and Fulman~\cite{BryanF} (see also~\cite{AbdesselamBDV} and references therein).
This is a holomorphic function of $z$ in the disk $|z|<1$, and we can extract coefficients using Cauchy's formula
\begin{equation}
\mathcal{H}_{\ell,n}(x)=\frac{1}{2i\pi}\oint_{C(r)}
z^{-n}\ \mathcal{G}_{\ell}(x,z)\ \frac{{\rm d}z}{z}\ ,
\label{Cauchyeq}
\end{equation}
with integration over the circle of radius $r\in(0,1)$ around the origin with counterclockwise orientation.
We then optimize the radius $r=e^{-t}$ or rather the value of the associated parameter $t\in(0,\infty)$, in order to minimize the supremum over the contour of the modulus of the integrand $|z^{-n}\mathcal{G}_{\ell}(x,z)|$.
We then split the coefficient of interest as a product
\begin{equation}
\mathcal{H}_{\ell,n}(x)=\mathcal{P}_{\ell,n}(x,t)\ \mathcal{J}_{\ell,n}(x,t)
\label{Hdecompeq}
\end{equation}
made of a prefactor
\begin{equation}
\mathcal{P}_{\ell,n}(x,t):=e^{nt}\mathcal{G}_{\ell}(x,e^{-t})
\label{Pdefeq}
\end{equation}
and an integral
\[
\mathcal{J}_{\ell,n}(x,t):=\int_{-\pi}^{\pi}j_{\ell,n}(x,t,\theta)\ \frac{{\rm d}\theta}{2\pi}
\]
where 
\[
j_{\ell,n}(x,t,\theta):=e^{-in\theta}\times\frac{\mathcal{G}_{\ell}(x,e^{-t+i\theta})}{\mathcal{G}_{\ell}(x,e^{-t})}\ .
\]
The proof of Proposition \ref{keyprop} relies on careful asymptotics of $\ln \mathcal{H}_{\ell,n}(x_n)$ for two different sequences $x_n$. One is the constant sequence $x_n=x$, and the other is the sequence $x_n=x e^{\frac{s}{b_n}}$. We use the exact contour radius optimizers $t_n$, as defined by the constant $x_n$ sequence, in order to also analyze the second non-constant sequence, where the $t_n$ are only approximate optimizers.
Since this may be a result of interest, in itself, we record the constant sequence asymptotics in the next proposition.
\begin{proposition}\label{logHprop}
For any $\ell\ge 2$, and any $x>0$, we have the asymptotic equivalence, as $n\rightarrow\infty$,
\[
\ln \mathcal{H}_{\ell,n}(x)\sim \frac{\ell}{\ell-1}\times(x\mathcal{K}_{\ell})^{\frac{1}{\ell}}\times n^{\frac{\ell-1}{\ell}}\ .
\]
\end{proposition}
When $x=1$, this immediately follows from the much more delicate asyptotics for
$ \mathcal{H}_{\ell,n}(1)$ itself, rather than its logarithm obtained in~\cite{BringmannFH}.
In order to control the integrals $\mathcal{J}_{\ell,n}$, we use a saddle point analysis with only one major arc corresponding to $|\theta|\le t_n$, and one minor arc corresponding to $|\theta|>t_n$. Finally, the prefactors $\mathcal{P}_{\ell,n}$ are controlled using leading $t\rightarrow 0^{+}$ asymptotics of the following remarkable multiple series.
For any $\ell\ge 0$, any complex numbers $\alpha_1,\ldots,\alpha_{\ell}$, and any $t\in(0,\infty)$, let
\[
Z_{\alpha_1,\ldots,\alpha_{\ell}}^{[\ell]}(t):=\sum_{\delta_1,\ldots,\delta_{\ell}=1}^{\infty}\delta_1^{\alpha_1-1}\cdots\delta_{\ell}^{\alpha_{\ell}-1} e^{-\delta_1\cdots\delta_{\ell}t}\ .
\]
By definition, $Z_{\varnothing}^{[0]}=e^{-t}$, while $Z_{\alpha}^{[1]}(t)=\sum_{\delta=1}^{\infty}\delta^{\alpha-1}e^{-\delta t}$ is a discrete analogue of $\int_{1}^{\infty}u^{\alpha-1}e^{-ut}{\rm d}u$ which is $t^{-\alpha}\Gamma(\alpha,t)$, in terms of the upper incomplete Gamma function. This explains the more convenient shift by 1 of the exponents for the $\delta$'s.

\subsection{Relation to other work and possible directions for further inquiry}
\label{subsec:relations}

The original motivation for the present article was the log-concavity conjecture by the first author for the numbers $A(\ell,n,k)$ with respect to $k$~\cite{Abdesselam} (see also~\cite[Challenge 3]{HeimN2} for $\ell=2$). This, of course, is equivalent to the log-concavity of the modified sequence $A(\ell,n,k)x^k$ which is proportional to the PMF of the random variables $\mathsf{K}_{\ell,n}$. Since the density of the standard Gaussian is the quintessential log-concave function, our CLT provides, in a weak sense, some confirmation for the conjecture. When $\ell\ge 3$, the only previously available results about the conjecture concern the ``dilute polymer gas'' regime where values of $k$ are very close to $n$, i.e., where nontrivial joint orbits are small and rare inside the environment $[n]$ (see~\cite[Cor. 4]{HeimN1},~\cite[Prop. 3.1]{AbdesselamBDV}, and~\cite{Tripathi}).
Note that for $\ell=2$, the log-concavity was established for a wide range of values of $k$ when $n$ is large enough ~\cite{Zhang}.

Obvious directions to explore beyond our CLT result are the study of moderate and large deviations for the $\mathsf{K}_{\ell,n}$, and upgrading the CLT to a local CLT. The latter would require asymptotics of the $A(\ell,n,k)$, rather than $\mathcal{H}_{\ell,n}(x)$, when both $n$ and $k$ are suitably large (see~\cite{StarrAbund,Abdesselam2025} for some results in this direction). Another interesting possible upgrade would be a functional CLT for a process associated to our $\mathbb{P}_{\ell,n,x}$ measures, in the spirit of~\cite{Hansen} (see~\cite{Pitman} for an in-depth study of such processes). 

Our proof of the CLT uses brute force asymptotic methods, and it would be desirable to find a more probabilistic proof, e.g., using a generalization of the Feller coupling idea. For $\ell=1$, $x\neq 1$, a generalization of this coupling for the Ewens measure was given in~\cite{ArratiaBT}. Within the bijective combinatorial approach to the Bryan-Fulman formula~\cite{AbdesselamBDV}, a detailed description of the joint orbits of a tuple $(\sigma_1,\ldots,\sigma_{\ell})$ as twisted discrete $\ell$-dimensional tori was introduced. This could be used to provide an analogue of the canonical cycle writing of a single permutation, which is the starting point of the Feller coupling idea, and should result in an expression for the random variables $\mathsf{K}_{\ell,n}$ as sums of $n$ Bernoulli random variables, albeit not independent ones. An analysis of correlations, e.g., if proven to decay fast enough, may allow one to reprove our CLT in a more probabilistically natural way.

Using the results of~\cite{AbdesselamBDV}, it is easy to generalize the Bryan-Fulman formula, and our measures $\mathbb{P}_{\ell,n,x}$ so the weight per joint orbit $x$ can depend on the size or other finer features of the orbit. In the case of a single permutation, the only (conjugation) isomorphism invariant of an orbit, i.e., cycle is the size. For $\ell\ge 2$, the discrete $\ell$-tori have more ``moduli'' such as directional linear sizes, and a parameter (denoted $z$ in~\cite{AbdesselamBDV}) for the twist used when gluing the ends of a cylinder based on an $(\ell-1)$-dimensional similar discrete torus. This should give rise to more general measures than $\mathbb{P}_{\ell,n,x}$. Allowing only a dependence of $x$ on the size of an orbit is straightforward, since the $A(\ell,n,k)$ can be written in terms of the $A(\ell,n,1)$ (see~\cite[Thm. 1.1]{AbdesselamBDV}), e.g., using Bell polynomials (see~\cite[Ch. 1]{Pitman}). This is also called a polymer gas representation in rigourous statistical mechanics~\cite{GruberK}. 

As emphasized in~\cite{Olshanski}, random permutations via the $n$-dimensional representation of $\mathfrak{S}_n$ by permutation matrices, can be seen as random matrix models which are discrete analogues of the circular unitary ensemble. The random tuples studied in this article can thus also be seen as discrete models of random commuting unitary matrices. The literature on random commuting matrices, especially with continuous models, is rather scarce (see however~\cite{McCarthy1,McCarthy2}).
One reason for this is that the commuting variety is a complicated mathematical object. One way to approach such models is to consider multi-matrix models with interaction terms given by squares of commutators, with a coupling constant $\beta$ which, when taken to infinity, forces the pairwise commutation. Such models were considered in physics~\cite{KazakovKN,KazakovZ}, and were also the object of recent rigorous mathematical work in~\cite{GuionnetMS}.
Note that counting $\ell$-tuples of commuting permutations is related to counting isomorphism classes under simultaneous conjugation of $(\ell-1)$-tuples of commuting permutations. Indeed, the number of such conjugacy classes is $\mathcal{H}_{\ell,n}(1)$, as follows from an easy application of Burnside's Lemma for the action of $\mathfrak{S}_n$ on $(\ell-1)$-tuples by simultaneous conjugation.
This shows for instance that $\mathcal{H}_{2,n}(1)=p(n)$ the number of integer partitions of $n$. 
Indeed the asymptotics in~\cite{BringmannFH} include, as a special case, the famous leading Hardy-Ramanujan asypmtotics
\[
\mathcal{H}_{2,n}(1)\sim
\frac{1}{4n\sqrt{3}}\exp\left(\pi\sqrt{\frac{2n}{3}}\right)\ ,
\]
recently revisited in~\cite{StarrHR}.
For $\ell\ge 3$, the Bringmann-Franke-Heim formulas are more delicate and rely on the deep generalization~\cite{BridgesBBF} of the theorem by Meinardus to the case of multiple poles for the relevant Dirichlet series. By contrast, this article is self-contained, and we tried to make it accessible even to readers who are not experts in the techniques of~\cite{deBruijn} and~\cite{FlajoletS}.
In light of the remark about conjugacy classes of shorter tuples, we note that if instead of random tuples of commuting matrices, one is interested in their isomorphism classes under simultaneous conjugation, such a model would live on the quotient of the commuting variety which is related to the Hilbert scheme of points in affine space (see the appendices of~\cite{Thomas} for a gentle introduction). Note that the matrix tuples featuring in the Hilbert scheme have the extra requirement that they must have a jointly cyclic vector. This is the linear algebra analogue of requiring $\kappa_{\ell,n}(\sigma_1,\ldots,\sigma_{\ell})=1$.

Note that a random tuple of commuting permutations $(\sigma_1,\ldots,\sigma_{\ell})$ is the same as a random group homomorphism $\varphi\in{\rm Hom}(\mathbb{Z}^{\ell},\mathfrak{S}_n)$. There is a vast body of literature studying, more generally, random homomorphisms $\varphi\in{\rm Hom}(\Gamma,\mathfrak{S}_n)$ where $\Gamma$ is a discrete finitely presented group (see, e.g.,~\cite{MageePvH} and references therein). There, the emphasis is on groups that are far from commutative, unlike $\mathbb{Z}^{\ell}$, but the type of questions investigated is quite different from the 
focus of this article. Most of the results concern the uniform measure on ${\rm Hom}(\Gamma,\mathfrak{S}_n)$. One can easily define analogues of our $\mathbb{P}_{\ell,n,x}$ in this setting by weighting a homomorphism $\varphi$ by $x^{\kappa_n(\varphi(\Gamma))}$. Whether these Ewens-like measures could be of use in the area of~\cite{MageePvH} remains to be seen. A class of groups $\Gamma$ which may be worth looking at, from the CLT perspective of our article, are right-angled Artin groups which interpolate between $\mathbb{Z}^{\ell}$ and the free group on $\ell$ elements (see, e.g.,~\cite{MageeT}).

Finally, we would like to mention that our work is very closely related to the work of Ercolani and Ueltschi~\cite{ErcolaniU} as well as that of Maples, Nikeghbali and Zeindler~\cite{MaplesNZ}. Indeed, there has been much recent activity (see~\cite{ElboimG} and references therein) studying the so-called generalized Ewens measures, i.e., probability distributions on single permutations $\sigma\in\mathfrak{S}_n$, with probability mass function
\[
\mathbb{P}(\sigma)=\frac{1}{h_n n!}\prod_{j=1}^{n}\theta_j^{c_j(\sigma)}\ ,
\]
where $c_j(\sigma)$ is the number of cycles of length $j$ in the permutation $\sigma$, and the normalization factor is
\[
h_n=\frac{1}{n!}\sum_{\sigma\in\mathfrak{S}_n}\prod_{j=1}^{n}\theta_j^{c_j(\sigma)}\ .
\]
The measure thus depends on the choice of weights $(\theta_j)_{j\ge 1}$. The basic Ewens measure corresponds to constant weights $\theta_j=x$. It is an easy exercise to show that our random variables $\mathsf{K}_{\ell,n}$ can be recovered as the number of cycles of the single random permutation $\sigma$ sampled according to such a generalized Ewens measure for the choice of weights
\[
\theta_j=\frac{x\ A(\ell,j,1)}{(j-1)!}\ .
\]
For $\ell\ge 2$, these are, on average, polynomially growing weights, i.e., belong to a less studied class of models in comparison to the case of logarithmically growing weights~\cite{ArratiaBTbook}. Nevertheless, there are some results on polynomially growing weights. For certain specific but different $\theta_j$ that on average grow like ours, mean asymptotics consistent with ours have been obtained in~\cite{ErcolaniU}, while a CLT for these simpler polynomially growing weights was derived in~\cite{MaplesNZ}.
Note that our weights grow like a power of $j$ only in an averaged or Cesaro sense. When $\ell=2$, we have $\theta_j=x \sigma(j)$ where $\sigma(j)=\sum_{q|j}q$ is the sum of divisors function of number theory. Simply showing sharp upper bounds on such a function is tantamount to the Riemann Hypothesis~\cite{Robin}. Note that it is possible that one may be able to prove our CLT using the general theorem in~\cite{MaplesNZ}, but this would require showing the generating function for our $\theta_j$ is log-admissible \`a la Hayman.
This typically is a nontrivial task (see, e.g.,~\cite[Theorem 6.2]{CantonFFM}). Although exploring a different direction of generalization, one should also mention the work~\cite{AhmadiGW} on asymptotics for generating series of similar flavor to the one considered by Bryan and Fulman.

\begin{remark}
It is well-known that to prove convergence in distribution, it suffices to prove convergence of the characteristic functions.
This is known as L\'evy's continuity criterion.
For example, this is the method used in references such as 
\cite{ElboimG,ErcolaniU,McCarthy2}.
In the present article we instead use the moment generating function.
In our context, the moment generating function exists, and we use monotonicity
for some of our bounds.
In order to use the characteristic function, one would need to extend these bounds.
That would be an excellent step towards establishing the methods in greater generality,
since some other types of examples might not have a convergent moment generating function on the entire
real line.
But we have not carried out this extension.
\end{remark}

\section{Preliminaries on the $Z$ functions}\label{prelimsec}

We begin with establishing some basic properties of the functions $Z_{\alpha_1,\ldots,\alpha_{\ell}}^{[\ell]}(t)$ which play an important role in this article. The case $Z_{\varnothing}^{[0]}$ is trivial, so let us assume $\ell\ge 1$. Let us also first focus on the case where the $\alpha$'s are real numbers. For any $\beta>0$ that also satisfies $\beta>\alpha_i$, for all $i\in[\ell]$, and for all $t>0$, we can write
\begin{eqnarray}
Z_{\alpha_1,\ldots,\alpha_{\ell}}^{[\ell]}(t) & = & t^{-\beta}
\sum_{\delta_1,\ldots,\delta_{\ell}=1}^{\infty}
\delta_1^{\alpha_1-1-\beta}\cdots\delta_{\ell}^{\alpha_{\ell}-1-\beta}\ (\delta_1\cdots\delta_{\ell}t)^{\beta} e^{-\delta_1\cdots\delta_{\ell}t} \nonumber \\
 & \le & t^{-\beta}\ \zeta(\beta+1-\alpha_1)\cdots \zeta(\beta+1-\alpha_{\ell})\ \beta^{\beta} e^{-\beta} \label{basicboundeq}\\
 & < & \infty\ , \nonumber
\end{eqnarray}
where we used the well know fact $\sup_{u>0}u^{\beta}e^{-u}=\beta^{\beta} e^{-\beta}$ because $\beta>0$.

Now allowing the $\alpha$'s to be complex, we have
\[
\sum_{\delta_1,\ldots,\delta_{\ell}=1}^{\infty}\left|
\delta_1^{\alpha_1-1}\cdots\delta_{\ell}^{\alpha_{\ell}-1} 
e^{-\delta_1\cdots\delta_{\ell}t}
\right|
=Z_{{\rm Re}(\alpha_1),\ldots,{\rm Re}(\alpha_{\ell})}^{[\ell]}(t)<\infty
\]
from the previous estimate, say with $\beta=1+\max(0,{\rm Re}\ \alpha_1,\ldots,{\rm Re}\ \alpha_{\ell})$.
Hence the series defining $Z_{\alpha_1,\ldots,\alpha_{\ell}}^{[\ell]}(t)$ converge absolutely and these are entire analytic functions of the $\alpha$'s in $\mathbb{C}^{\ell}$. However, we will be focusing on the dependence
on $t\in(0,\infty)$, with the $\alpha$'s fixed.
Since for any $T>0$, and any $t\in[T,\infty)$,
\[
\left|-\delta_1^{\alpha_1}\cdots\delta_{\ell}^{\alpha_{\ell}} e^{-\delta_1\cdots\delta_{\ell}t}
\right|\le
\delta_1^{{\rm Re}(\alpha_1)}\cdots\delta_{\ell}^{{\rm Re}(\alpha_{\ell})} e^{-\delta_1\cdots\delta_{\ell}T}
\]
and $Z_{{\rm Re}(\alpha_1)+1,\ldots,{\rm Re}(\alpha_{\ell})+1}^{[\ell]}(T)<\infty$, the corollary of the dominated convergence theorem pertaining to differentiation under the integral sign applies. Hence, $Z_{\alpha_1,\ldots,\alpha_{\ell}}^{[\ell]}(t)$ is differentiable, and for all $t\in(0,\infty)$, we have
\begin{equation}
\frac{{\rm d}}{{\rm d}t}Z_{\alpha_1,\ldots,\alpha_{\ell}}^{[\ell]}(t)=-Z_{\alpha_1+1,\ldots,\alpha_{\ell}+1}^{[\ell]}(t)\ .
\label{Zderivativeeq}
\end{equation}
By iteration, we see that $Z_{\alpha_1,\ldots,\alpha_{\ell}}^{[\ell]}(t)$ is $C^{\infty}$ on $(0,\infty)$.

By splitting the factor $e^{-\delta_1\cdots\delta_{\ell}t}$ in two, and using $\delta_1\cdots\delta_{\ell}\ge 1$ because the $\delta$'s are $\ge 1$, we immediately obtain the inequality
\begin{equation}
\left|Z_{\alpha_1,\ldots,\alpha_{\ell}}^{[\ell]}(t)\right|\le e^{-\frac{t}{2}}
\ Z_{{\rm Re} (\alpha_1),\ldots,{\rm Re} (\alpha_{\ell})}^{[\ell]}
\left(\frac{t}{2}\right)\ .
\label{freedecayeq}
\end{equation}

\begin{lemma}\label{Zboundlem}
For any $\ell\ge 0$, any $\alpha_1,\ldots,\alpha_{\ell}$ in $\mathbb{C}$, and any $\beta>\max(0, {\rm Re}\ \alpha_1,\ldots,{\rm Re}\ \alpha_{\ell})$, there exists $c_1,c_2>0$, such that for all $t>0$, we have
\[
\left|Z_{\alpha_1,\ldots,\alpha_{\ell}}^{[\ell]}(t)\right|\le c_1 t^{-\beta} e^{-c_2 t}\ .
\]
\end{lemma}

\noindent{\bf Proof:}
The result is trivial for $\ell=0$, with $c_1=\frac{1}{2}$ and $c_2=(2\beta)^{\beta} e^{-\beta}$, again using $\sup_{u>0}u^{\beta}e^{-u}=\beta^{\beta} e^{-\beta}$.
For $\ell\ge 1$, we first use (\ref{freedecayeq}), and then we bound 
$Z_{{\rm Re} (\alpha_1),\ldots,{\rm Re} (\alpha_{\ell})}^{[\ell]}
\left(\frac{t}{2}\right)$ using the estimate (\ref{basicboundeq}) for $\frac{t}{2}$ instead of $t$. The stated inequality then holds with $c_1=\frac{1}{2}$ and
\[
c_2=\zeta(\beta+1-\alpha_1)\cdots \zeta(\beta+1-\alpha_{\ell})\ (2\beta)^{\beta} e^{-\beta}\ ,
\]
as wanted.
\qed

In order to obtain the leading asymptotics for the $Z_{\alpha_1,\ldots,\alpha_{\ell}}^{[\ell]}(t)$ when $t\rightarrow 0^{+}$, we need the following elementary lemma which controls the Euler--Maclaurin approximation of integrals on $(0,\infty)$ which are possibly improper near zero.

\begin{lemma}\label{EulerMaclem}
Let $f\colon(0,\infty)\rightarrow\mathbb{C}$ be a differentiable function. Suppose there exist $\gamma>0$, such that
$f(t) \ll t^{\gamma-1}$ as $t \to 0^+$ and $f(t)\ll t^{-\gamma-1}$ as $t \to \infty$.
Then $f$ is integrable on $(0,\infty)$, the series $\sum_{\delta=1}^{\infty}f(\delta t)$ converges absolutely for all $t>0$, and we have
\[
\lim\limits_{t\rightarrow 0^{+}}\ t\sum_{\delta=1}^{\infty}f(\delta t)
=\int_{0}^{\infty}f(u)\ {\rm d}u\ .
\]
\end{lemma}

\noindent{\bf Proof:}
In P\'olya and Szeg\"o \cite{PolyaSzego}, Problem 159 on page 87, it is established that the limit
$\lim_{N \to \infty} \frac{1}{N} \sum_{n=1}^{\infty} f(n/N)$ converges to $\int_0^{\infty} f(u)\ {\rm d}u$
as long as the conditions are satisfied. Their problem is somewhat more general than this: we have used it for the regular
sequence being $r_n=n$ which has convergence exponent $\lambda=1$.
(We thank an anonymous referee for pointing out this reference to us.)
\qed

Note that the functions $Z_{\alpha_1,\ldots,\alpha_{\ell}}^{[\ell]}(t)$ are symmetric in the $\alpha$'s and there is no harm in restricting to the case where the real parts are decreasing, i.e.,
${\rm Re}\ \alpha_1\ge\cdots\ge {\rm Re}\ \alpha_{\ell}$.
In general, when the inequalities are not strict, the $t\rightarrow 0^{+}$ asymptotics contain logarithms. The next key lemma provides leading asymptotics, without logarithms,
in sufficient generality for the needs of this article.

\begin{lemma}\label{Zasymlem}
Let $\ell\ge 1$, and suppose the complex numbers $\alpha_1,\ldots,\alpha_{\ell}$ are such that ${\rm Re}\ \alpha_1>{\rm Re}\ \alpha_2\ge\cdots\ge {\rm Re}\ \alpha_{\ell}$ and ${\rm Re}\ \alpha_1>0$. Then, as $t\rightarrow 0^{+}$, we have the asymptotic equivalence
\[
Z_{\alpha_1,\ldots,\alpha_{\ell}}^{[\ell]}(t)\sim
\Gamma(\alpha_1)\zeta(\alpha_1+1-\alpha_2)\zeta(\alpha_1+1-\alpha_3)\cdots\zeta(\alpha_1+1-\alpha_{\ell})\ t^{-\alpha_1}\ ,
\]
in terms of the Euler gamma and the Riemann zeta functions.
Note that when $\ell=1$ the only hypothesis is ${\rm Re}\ \alpha_1>0$, and there are no zeta factors in the conclusion.
\end{lemma}

\noindent{\bf Proof:}
The sums defining the $Z$ functions having been shown to be absolutely convergent, we can use the discrete Fubini theorem to write
\[
Z_{\alpha_1,\ldots,\alpha_{\ell}}^{[\ell]}(t)=\sum_{\delta_1=1}^{\infty}
\delta_1^{\alpha_1-1}Z_{\alpha_2,\ldots,\alpha_{\ell}}^{[\ell-1]}(\delta_{1}t)
=t^{-\alpha_1}\times t\sum_{\delta_1=1}^{\infty}f(\delta_1 t)\ ,
\]
with the new function on $(0,\infty)$ given by
\[
f(t):=t^{\alpha_1-1} Z_{\alpha_2,\ldots,\alpha_{\ell}}^{[\ell-1]}(t)\ .
\]
The functions $f$ 
will be shown to satisfy the required bounds for the application of Lemma \ref{EulerMaclem}.
By hypothesis, ${\rm Re}\ \alpha_1>\max(0,{\rm Re}\ \alpha_2)$, and it is therefore possible to choose some real number $\beta$ such that
${\rm Re}\ \alpha_1>\beta>\max(0,{\rm Re}\ \alpha_2)$.
We now use Lemma \ref{Zboundlem} which provides constants $c_{1,1},c_{1,2}>0$ such that, for all $t>0$,
\[
\left|Z_{\alpha_2,\ldots,\alpha_{\ell}}^{[\ell-1]}(t)\right|
\le c_{1,1} t^{-\beta} e^{-c_{1,2}t}\ .
\]
We then easily see that the estimates needed for the use of Lemma \ref{EulerMaclem} hold with $\gamma:={\rm Re}(\alpha_1)-\beta>0$.
Hence, Lemma \ref{EulerMaclem} shows that
\[
\lim\limits_{t\rightarrow 0^{+}} t^{\alpha_1}Z_{\alpha_1,\ldots,\alpha_{\ell}}^{[\ell]}(t)=\int_{0}^{\infty}f(u)\ {\rm d}u\ .
\]
All that remains is to compute the last integral.
We first consider the case where all the $\alpha$'s are real, so one can commute sums and integrals via Tonelli's theorem, without fear. 
We then have
\begin{eqnarray}
\int_{0}^{\infty}f(u)\ {\rm d}u & = & \int_{0}^{\infty}
u^{\alpha_1-1}\left(
\sum_{\delta_2,\ldots,\delta_{\ell}=1}^{\infty}\delta_2^{\alpha_2-1}\cdots\delta_{\ell}^{\alpha_{\ell}-1} e^{-\delta_2\cdots\delta_{\ell}u}
\right)\ {\rm d}u \nonumber \\
 &= & \sum_{\delta_2,\ldots,\delta_{\ell}=1}^{\infty}\delta_2^{\alpha_2-1}\cdots\delta_{\ell}^{\alpha_{\ell}-1}
\int_{0}^{\infty}u^{\alpha_1-1}e^{-\delta_2\cdots\delta_{\ell}u}\ {\rm d}u
\nonumber \\
 & = & \sum_{\delta_2,\ldots,\delta_{\ell}=1}^{\infty}\delta_2^{\alpha_2-1}\cdots\delta_{\ell}^{\alpha_{\ell}-1}\times\frac{\Gamma(\alpha_1)}{(\delta_2\cdots\delta_{\ell})^{\alpha_1}}
\nonumber \\
 &= & \Gamma(\alpha_1)\sum_{\delta_2,\ldots,\delta_{\ell}=1}^{\infty}
\frac{1}{\delta_2^{\alpha_1+1-\alpha_2}}\times\cdots\times\frac{1}{\delta_{\ell}^{\alpha_1+1-\alpha_{\ell}}} \nonumber \\
 &= & \Gamma(\alpha_1)\zeta(\alpha_1+1-\alpha_2)\zeta(\alpha_1+1-\alpha_3)\cdots\zeta(\alpha_1+1-\alpha_{\ell}) \label{intcompeq}
\end{eqnarray}
which is finite.
We then move on to the case where the $\alpha$'s are allowed to be complex. The above computation with the real parts of the $\alpha$'s instead of the $\alpha$'s themselves
shows that all commutations of summations and integrations are justified when redoing the same computation in the complex case. Therefore the last formula (\ref{intcompeq}) holds for the full scope of the statement of the present lemma, which is now established.
\qed

In this article, we will only need the above result for some descending staircase shapes. It is convenient, for $\ell\ge 1$ and $m\in\mathbb{Z}$, to introduce the abbreviated notation
\[
Z_{m}^{[\ell]}(t):=Z_{m,m-1,\ldots,m-\ell+1}^{[\ell]}(t)\ .
\]
For ease of reference, we record below the particular asymptotics needed for the remainder of this article. Thanks to Lemma \ref{Zasymlem}, when $\ell\ge 2$, the following $t\rightarrow 0^{+}$ asymptotic equivalences hold
\begin{eqnarray}
Z_{\ell-1}^{[\ell]}(t) & \sim & \frac{\mathcal{K}_{\ell}}{\ell-1}\ t^{-\ell+1}
\label{Zminus1asymeq} \\
Z_{\ell}^{[\ell]}(t) & \sim & \mathcal{K}_{\ell}\ t^{-\ell}
\label{Zasymeq} \\
Z_{\ell+1}^{[\ell]}(t) & \sim & \ell\mathcal{K}_{\ell}\ t^{-\ell-1}
\label{Zplus1asymeq} \\
Z_{\ell+2}^{[\ell]}(t) & \sim & \ell(\ell+1)\mathcal{K}_{\ell}\ t^{-\ell-2} \ .
\label{Zplus2asymeq}
\end{eqnarray}
The constant $\mathcal{K}_{\ell}$ is the one defined in (\ref{Kdefeq}). 
The second function in (\ref{Zasymeq}) deserves special attention, which explains which of the four constants we decided to give a name to.

From the previous facts established in this section, it is clear that $Z_{\ell}^{[\ell]}(t)$ is a smooth strictly decreasing function from $(0,\infty)$ to itself. The above asymptotics imply $\lim_{t\rightarrow 0^{+}}Z_{\ell}^{[\ell]}(t)=\infty$. Lemma \ref{Zboundlem}, say with $\beta=\ell+1$, provides exponential decay at infinity, and implies $\lim_{t\rightarrow\infty}Z_{\ell}^{[\ell]}(t)=0$.
Therefore, $Z_{\ell}^{[\ell]}(t)$ is a bijection from $(0,\infty)$ onto itself, and the inverse function $(Z_{\ell}^{[\ell]})^{-1}:(0,\infty)\rightarrow(0,\infty)$ is well defined and bijective too.
The last ingredient needed from this section is the following calculus exercise.

\begin{lemma}\label{calclem}
As $u\rightarrow\infty$, we have the asymptotic equivalence
\[
(Z_{\ell}^{[\ell]})^{-1}(u)\sim \left(\frac{u}{\mathcal{K}_{\ell}}\right)^{-\frac{1}{\ell}}\ .
\]
\end{lemma}
 
\noindent{\bf Proof:}
We have
\[
\left(\frac{u}{\mathcal{K}_{\ell}}\right)^{\frac{1}{\ell}}
(Z_{\ell}^{[\ell]})^{-1}(u)=\left[\frac{Z_{\ell}^{[\ell]}(t)}{\mathcal{K}_{\ell}t^{-\ell}}\right]^{\frac{1}{\ell}}
\]
with $t=(Z_{\ell}^{[\ell]})^{-1}(u)$. Since $t\rightarrow 0^{+}$ when $u\rightarrow\infty$, and since $\lim_{t\rightarrow 0^{+}}\mathcal{K}_{\ell}^{-1}t^{\ell}Z_{\ell}^{[\ell]}(t)=1$, the left-hand side of the above equation goes to $1$, and the lemma follows.
\qed

\section{Preparation for the saddle point analysis}

For $u\in\mathbb{C}\backslash (-\infty,0]$ we will use ${\rm Log}\ u:=\ln |u|+i{\rm Arg}\ u$ for the principal branch of the complex logarithm, namely, with $-\pi<{\rm Arg}\ u<\pi$. We of course have 
\begin{equation}
{\rm Log}(1-u)=-\sum_{\delta=1}^{\infty}\frac{u^{\delta}}{\delta}\ ,
\label{Logeq}
\end{equation}
when $|u|<1$.
For $\ell\ge 2$, and $z$ in the open unit disk around the origin, we define
\begin{eqnarray}
\mathcal{L}_{\ell}(z) & := & -\sum_{\delta_1,\ldots,\delta_{\ell-1}=1}^{\infty}
\delta_{1}^{\ell-2}\delta_{2}^{\ell-3}\cdots\delta_{\ell-2}\ {\rm Log}(1-z^{\delta_1\cdots\delta_{\ell-1}}) \label{Ldefeq} \\
 & = & \sum_{\delta_1,\ldots,\delta_{\ell}=1}^{\infty}
\delta_{1}^{\ell-2}\delta_{2}^{\ell-3}\cdots\delta_{\ell-2}
\delta_{\ell-1}^{0}\delta_{\ell}^{-1}
z^{\delta_1\cdots\delta_{\ell}} \label{Lniceeq}\ ,
\end{eqnarray}
where we used (\ref{Logeq}) to introduce an extra last summation index. The use of the discrete Fubini theorem is justified because the series (\ref{Lniceeq}) converges absolutely, since $Z_{\ell-1}^{[\ell]}(-\ln |z|)<\infty$, by the results of \S\ref{prelimsec}.
Therefore, the series in (\ref{Ldefeq}) converges absolutely, which easily implies
\[
\sum_{\delta_1,\ldots,\delta_{\ell-1}=1}^{\infty}
\left|
\left(1-z^{\delta_1\cdots\delta_{\ell-1}}\right)^{-x\delta_1^{\ell-2}\delta_2^{\ell-3}\cdots\delta_{\ell-2}}-1
\right|<\infty\ .
\]
In other words, the product on the right-hand side of (\ref{BryanFeq}) is absolutely convergent. Hence, this right-hand side is a holomorphic function of $z$, for $|z|<1$,
and the formal power series in $\mathbb{C}[[z]]$ corresponding to its Taylor expansion at the origin is the infinite product, in the sense of formal power series, of the Taylor series at the origin of the individual factors. This product was shown by Bryan and Fulman~\cite{BryanF} to be the formal power series given by the left-hand side of (\ref{BryanFeq}). This therefore shows that the left-hand side of (\ref{BryanFeq}) also converges in the disk $|z|<1$ and is equal to the infinite product on the right-hand side, not just as formal power series, but as holomorphic functions.

Now that (\ref{Cauchyeq}) is justified with $r=e^{-t}$ for some $t>0$ to be chosen shortly, we look at the modulus 
of the integrand of equation (\ref{Cauchyeq}), namely
$|z^{-n}\mathcal{G}_{\ell}(x,z)|=\exp(nt+x{\rm Re}\ \mathcal{L}_{\ell}(z))$. For $z=e^{-t+i\theta}$, with $-\pi<\theta<\pi$, we have
\[
{\rm Re}\ \mathcal{L}_{\ell}(z)=
\sum_{\delta_1,\ldots,\delta_{\ell}=1}^{\infty}
\delta_{1}^{\ell-2}\delta_{2}^{\ell-3}\cdots\delta_{\ell-2}
\delta_{\ell-1}^{0}\delta_{\ell}^{-1}
\ e^{-\delta_1\cdots\delta_{\ell}t}\ \cos(\delta_1\cdots\delta_{\ell}\theta)\ ,
\]
which is clearly maximal when $\theta=0$. The maximum over the contour of the modulus of the integrand is thus given by $\exp(nt+xZ_{\ell-1}^{[\ell]}(t))$. Taking the derivative of the logarithm of the last expression, we see that it is minimal when the equation
\[
n-xZ_{\ell}^{[\ell]}(t)=0
\]
holds, i.e., when
$t=(Z_{\ell}^{[\ell]})^{-1}\left(\frac{n}{x}\right)$.
We go ahead and use this to define the sequence
\begin{equation}
t_n:=(Z_{\ell}^{[\ell]})^{-1}\left(\frac{n}{x}\right)\ .
\label{tndefeq}
\end{equation}
By Lemma \ref{calclem}, when $n\rightarrow\infty$, we have
\begin{equation}
t_n\sim \left(\frac{n}{x\mathcal{K}_{\ell}}\right)^{-\frac{1}{\ell}}\ .
\label{tnasymeq}
\end{equation}
Recall that for the needs of Proposition \ref{keyprop}, we are not only considering a fixed $x$, but rather a sequence $x_n$.
Let $s\in\mathbb{R}$ be fixed, and define the sequence 
\begin{equation}
x_n:=x\ e^{\frac{s}{b_n}}\ ,
\label{xndefeq}
\end{equation}
where $(b_n)$ is the one in the statement of Proposition \ref{keyprop}.
Still in the setting of Proposition \ref{keyprop}, we introduce the notation
\begin{eqnarray}
\Psi_{\ell,n}(s) & := & \ln
\mathbb{E}\left[
\exp\left(s\left(\frac{\mathsf{K}_{\ell,n}-a_n}{b_n}\right)\right)
\right] \\
 & = & -\frac{sa_n}{b_n}
+\ln \mathcal{P}_{\ell,n}(x_n,t_n)
-\ln \mathcal{P}_{\ell,n}(x,t_n)
+R_{\ell,n}(s)\ ,
\label{Psiconvenienteq}
\end{eqnarray}
where the remainder is
\begin{equation}
R_{\ell,n}(s):= \ln\left[\frac{\mathcal{J}_{\ell,n}(x_n,t_n)}{\mathcal{J}_{\ell,n}(x,t_n)}\right]\ ,
\label{remaindereq}
\end{equation}
and where we used the multiplicative decomposition (\ref{Hdecompeq}).
Note that the $\mathcal{J}_{\ell,n}$ are positive real numbers given by ratios $\mathcal{H}_{\ell,n}/\mathcal{P}_{\ell,n}$, even if they are expressed as integrals of complex-valued functions.
We now rewrite the $\mathcal{J}_{\ell,n}$ integrals as
\begin{equation}
\mathcal{J}_{\ell,n}(x_n,t_n)=\int_{-\pi}^{\pi}e^{-q_{\ell,n}(x_n,t_n,\theta)}\ \frac{{\rm d}\theta}{2\pi}\ ,
\label{newJformulaeq}
\end{equation}
where, for $y,u>0$ and $-\pi<\theta<\pi$, the $q$ functions are defined by
\begin{eqnarray}
q_{\ell,n}(y,u,\theta) &:=&
in\theta-y\left(\mathcal{L}_{\ell}(e^{-u+i\theta})
-\mathcal{L}_{\ell}(e^{-u})\right) \nonumber \\
 &=& in\theta-y
\sum_{\delta_1,\ldots,\delta_{\ell}=1}^{\infty}
\delta_{1}^{\ell-2}\delta_{2}^{\ell-3}\cdots\delta_{\ell-2}
\delta_{\ell-1}^{0}\delta_{\ell}^{-1}
e^{-\delta_1\cdots\delta_{\ell}u}\left(e^{i\delta_1\cdots\delta_{\ell}\theta}-1\right)\ .
\label{qconvenienteq}
\end{eqnarray}
Using, for $\beta$ real, the identity ${\rm Re}\left(1-e^{i\beta}\right)=1-\cos(\beta)=2\sin^2\left(\frac{\beta}{2}\right)$, we get
\[
{\rm Re}\ q_{\ell,n}(y,u,\theta)=2y\sum_{\delta_1,\ldots,\delta_{\ell}=1}^{\infty}
\delta_{1}^{\ell-2}\delta_{2}^{\ell-3}\cdots\delta_{\ell-2}
\delta_{\ell-1}^{0}\delta_{\ell}^{-1}
\ e^{-\delta_1\cdots\delta_{\ell}u}
\ \sin^2\left(\frac{\delta_1\cdots\delta_{\ell}\theta}{2}\right)\ .
\]
Since all the terms are nonnegative, this gives the lower bound
\begin{equation}
{\rm Re}\ q_{\ell,n}(y,u,\theta)\ge 2y
\sum_{k=1}^{\infty}k^{\ell-2}\ e^{-ku}\ \sin^2\left(\frac{k\theta}{2}\right)
\label{qlowerbdeq}
\end{equation}
obtained by only keeping tuples of the form $(\delta_1,\ldots,\delta_{\ell})=(k,1,\ldots,1)$, and which is the main tool we use for the estimates of the next section.
We will split the integration domain for $\theta$ into one major arc region
\[
\mathcal{R}_{\rm maj}:=\{\theta\in (-\pi,\pi)\ |\ |\theta|\le t_n\}\ ,
\]
and one minor arc region
\[
\mathcal{R}_{\rm min}:=\{\theta\in (-\pi,\pi)\ |\ |\theta|> t_n\}\ .
\]
Although both based on (\ref{qlowerbdeq}), the estimates for the two regions will be treated differently.

\section{The saddle point analysis}

\subsection{Major arc estimates}

We use (\ref{qlowerbdeq}) to deduce, for $\theta\in\mathcal{R}_{\rm maj}$,
\begin{equation}
{\rm Re}\ q_{\ell,n}(x_n,t_n,\theta)\ge
2x_n
\sum_{k=1}^{N_n}k^{\ell-2}\ e^{-kt_n}\ \sin^2\left(\frac{k|\theta|}{2}\right)
\label{majasumeq}
\end{equation}
where we introduced the cutoff
\[
N_n:=\left\lfloor\frac{\pi}{t_n}\right\rfloor\ .
\]
This ensures that in the sum (\ref{majasumeq}) we have
\[
\frac{k|\theta|}{2}\le \frac{t_n}{2}\times \left\lfloor\frac{\pi}{t_n}\right\rfloor\le\frac{\pi}{2}\ .
\]
Since the sine function is convex on $\left[0,\frac{\pi}{2}\right]$, it is bounded below by the chord defined by the endpoints, i.e., we have $\sin u\ge \frac{2u}{\pi}$, for $0\le u\le\frac{\pi}{2}$. Note that the cutoff $N_n$ also ensures that in (\ref{majasumeq}) we have $t_nk\le \pi$.
Hence, we can write, for $n$ large enough,
\begin{eqnarray*}
{\rm Re}\ q_{\ell,n}(x_n,t_n,\theta) & \ge &
2x_n
\sum_{k=1}^{N_n}k^{\ell-2}\ e^{-\pi}\ \left(\frac{k|\theta|}{\pi}\right)^2 \\
 & \ge & \frac{2x_n\theta^2}{\pi^2e^{\pi}}\times\sum_{k=1}^{N_n}k^{\ell} \\
 & \ge & \frac{2x_n\theta^2}{\pi^2e^{\pi}}\times\int_{0}^{N_n} u^{\ell}\ {\rm d}u \\
 & \ge &  \frac{2x_n\theta^2 N_n^{\ell+1}}{\pi^2e^{\pi}(\ell+1)}\ .
\end{eqnarray*}
Since $x_n$ converges to $x$, and
\[
N_n=\left\lfloor\frac{\pi}{t_n}\right\rfloor
\sim\frac{\pi}{t_n}\sim\pi\times\left(\frac{n}{x\mathcal{K}_{\ell}}\right)^{\frac{1}{\ell}} \ ,
\]
there exists a constant $\eta_{\rm maj}>0$
 (possibly dependent on $\ell,x$) such that, for all $n$ and all $\theta\in\mathcal{R}_{\rm maj}$, we have the lower bound
\begin{equation}
{\rm Re}\ q_{\ell,n}(x_n,t_n,\theta)\ge
\eta_{\rm maj}\ \theta^{2}\ n^{\frac{\ell+1}{\ell}}\ .
\label{majalowerbdeq}
\end{equation}

\subsection{Minor arc estimates}

We return to (\ref{qlowerbdeq}) with generic arguments $y,u$ instead of the designated sequence terms $x_n,t_n$.
Since $\ell-2\ge 0$, we have $k^{\ell-2}\ge 1$, and therefore
\[
{\rm Re}\ q_{\ell,n}(y,u,\theta)\ge
2y\sum_{k=1}^{\infty}\ e^{-ku}\ \sin^2\left(\frac{k\theta}{2}\right)\ ,
\]
namely, ${\rm Re}\ q_{\ell,n}(y,u,\theta)\ge y\ {\rm Re}\  \rho(u,\theta)$ with
\[
\rho(u,\theta):=\sum_{k=1}^{\infty}e^{-ku}\left(1-e^{ik\theta}\right)\ ,
\] 
after again using the identity ${\rm Re}\left(1-e^{i\beta}\right)=2\sin^2\left(\frac{\beta}{2}\right)$, but the other way. With a little bit of algebra, we see that
\begin{eqnarray*}
\rho(u,\theta) &= & \frac{e^{-u}}{1-e^{-u}}-\frac{e^{-u+i\theta}}{1-e^{-u+i\theta}} \\
 &= &  \frac{1}{e^{u}-1}-\frac{e^{i\theta}}{e^u-e^{i\theta}} \\
 & = & \frac{e^u(1-e^{i\theta})}{(e^{u}-1)(e^u-e^{i\theta})} \\
 & = & \frac{e^u(e^u-e^{-i\theta}-e^u e^{i\theta}+1)}{(e^{u}-1)(e^{2u}-2e^u\cos \theta+1)}\ ,
\end{eqnarray*}
after multiplying, above and below, by the conjugate expression $e^u-e^{-i\theta}$.
Taking the real part and factoring the new numerator, we get
\[
{\rm Re}\ \rho(u,\theta) =
\frac{e^u(e^u+1)}{e^u-1}\times
\frac{1-\cos\theta}{e^{2u}-2e^u\cos \theta+1}\ .
\]
Note that for any fixed $c>1$, the function $v\mapsto\frac{1-v}{c^2-2cv+1}$ 
on the interval $(-\infty,\frac{c^2+1}{2c})$ which contains $[-1,1]$, has derivative 
$-\left(\frac{c-1}{c^2-2cv+1}\right)^2<0$.
Therefore, ${\rm Re}\ \rho(u,\theta)$ is a decreasing function of $\cos \theta$ or rather $\cos |\theta|$, by parity.
We now apply this to $y=x_n$, $u=t_n$, and we suppose $\theta\in\mathcal{R}_{\rm min}$, i.e., $\cos |\theta|\le \cos t_n$, with $n$ large enough so that $t_n<\pi$.
The last decreasing property then provides us with the lower bound
\[
{\rm Re}\ q_{\ell,n}(x_n,t_n,\theta)\ge
\frac{e^{t_n}(e^{t_n}+1)}{e^{t_n}-1}\times
\frac{1-\cos(t_n)}{e^{2t_n}-2e^{t_n}\cos(t_n)+1}\ .
\]
Since $t_n\rightarrow 0$, we have by expanding explicitly to second order
\begin{eqnarray*}
e^{2t_n}-2e^{t_n}\cos(t_n)+1 &= &
(1+2t_n+2t_n^2)-2\left(1+t_n+\frac{t_n^2}{2}\right)\left(1-\frac{t_n^2}{2}\right)+1+O(t_n^3)\ \\
 & = & 2t_n^2 +O(t_n^3)\ ,
\end{eqnarray*}
and therefore
\[
\frac{e^{t_n}(e^{t_n}+1)}{e^{t_n}-1}\times
\frac{1-\cos(t_n)}{e^{2t_n}-2e^{t_n}\cos(t_n)+1}
\sim \frac{2}{t_n}\times\frac{\left(\frac{t_n^2}{2}\right)}{2t_n^2}=\frac{1}{2t_n}\sim \frac{1}{2}\times\left(\frac{n}{x\mathcal{K}_{\ell}}\right)^{\frac{1}{\ell}}.
\]
As a result, there exists a constant $\eta_{\rm min}>0$
(possibly depending on $\ell,x$) such that for all $n$ and $\theta\in\mathcal{R}_{\rm min}$,
\[
{\rm Re}\ q_{\ell,n}(x_n,t_n,\theta)\ge \eta_{\rm min}\ n^{\frac{1}{\ell}}\ .
\]
This immediately implies the fractional exponential decay estimate
\begin{equation}
\left|
\int_{\mathcal{R}_{\rm min}}e^{-q_{\ell,n}(x_n,t_n,\theta)}\ \frac{{\rm d}\theta}{2\pi}\right|
\le \exp\left(-\eta_{\rm min}\ n^{\frac{1}{\ell}}\right)\ .
\label{minacrusheq}
\end{equation}

\subsection{Putting the estimates together and deriving the $\mathcal{J}$ integral asymptotics}

With $t_n$ and $x_n$ as in (\ref{tndefeq}) and (\ref{xndefeq}), define the new sequence
\begin{equation}
\lambda_n:=(x_n Z_{\ell+1}^{[\ell]}(t_n))^{-\frac{1}{2}}\ .
\label{lambdadefeq}
\end{equation}
From (\ref{xndefeq}), (\ref{tnasymeq}), and (\ref{Zplus1asymeq}), we easily get the asymptotic equivalent
\begin{equation}
\lambda_n\sim \frac{(x\mathcal{K}_{\ell})^{\frac{1}{2\ell}}}{\sqrt{\ell}}\times n^{-\left(\frac{\ell+1}{2\ell}\right)}\ .
\label{lambdaasymeq}
\end{equation}
We do the change of variable $\theta=\lambda_n\Theta$ in the portion of the integral $\mathcal{J}_{\ell,n}(x_n,t_n)$ over the major arc region $\mathcal{R}_{\rm maj}$, which results in
\[
\int_{\mathcal{R}_{\rm maj}}e^{-q_{\ell,n}(x_n,t_n,\theta)}\ \frac{{\rm d}\theta}{2\pi}=\frac{\lambda_n}{2\pi}\int_{\mathbb{R}}f_n(\Theta)\ {\rm d}\Theta
\]
involving the sequence of functions $f_n\colon\mathbb{R}\rightarrow\mathbb{C}$ define by
\[
f_n(\Theta):=\bbone\{|\Theta|\le t_n\lambda_n^{-1}\}
\ e^{-q_{\ell,n}(x_n,t_n,\lambda_n\Theta)}\ ,
\]
with $n$ assumed large enough so that $t_n<\pi$.
The notation $\bbone\{\cdots\}$ stands for the indicator function of the enclosed condition. 
We also note that $t_n\lambda_{n}^{-1}$, of order $n^{\frac{\ell-1}{2\ell}}$, goes to infinity.
From (\ref{majalowerbdeq}), we see that for all $\Theta\in\mathbb{R}$,
\[
|f_n(\Theta)|\le \exp\left(-\eta_{\rm maj}\ \lambda_n^2\ n^{\frac{\ell+1}{\ell}}\ \Theta^2\right)\ .
\]
Pick some constant $c$ such that
\[
0<c<\eta_{\rm maj}\times\frac{(x\mathcal{K}_{\ell})^{\frac{1}{\ell}}}{\ell}\ .
\]
By the asymptotics (\ref{lambdaasymeq}), for $n$ large enough we have
\[
\eta_{\rm maj}\ \lambda_n^2\ n^{\frac{\ell+1}{\ell}}\ge c\ ,
\]
and therefore the domination
\[
|f_n(\Theta)|\le e^{-c\Theta^2}
\]
by a fixed function which is integrable on $\mathbb{R}$.

We now compute the pointwise limits $\lim_{n\rightarrow\infty} f_n(\Theta)$.
For any $w\in\mathbb{C}$, we have by the Taylor formula with integral remainder
\[
e^w=1+w+\frac{w^2}{2}+R(w)
\]
with
\[
R(w)=\int_0^1 \frac{(1-\beta)^2}{2}\ w^3\ e^{\beta w} \ {\rm d}\beta
\]
which imples $|R(w)|\le \frac{|w|^3}{6}$ when $w$ is pure imaginary.
We insert this, with $w=i\lambda_n\delta_1\cdots\delta_{\ell}\Theta$, in the formula
(\ref{qconvenienteq}), which gives
\[
q_{\ell,n}(x_n,t_n,\lambda_n\Theta)=
i\lambda_n\Theta\left[n-x_n Z_{\ell}^{[\ell]}(t_n)\right]
+\frac{\Theta^2}{2}\times \lambda_n^2 x_n Z_{\ell+1}^{[\ell]}(t_n) 
+{\rm Err}\ ,
\]
with the explicit error term
\[
{\rm Err}:=-x_n\sum_{\delta_1,\ldots,\delta_{\ell}=1}^{\infty}
\delta_{1}^{\ell-2}\delta_{2}^{\ell-3}\cdots\delta_{\ell-2}
\delta_{\ell-1}^{0}\delta_{\ell}^{-1}
\ e^{-\delta_1\cdots\delta_{\ell}t_n}\ R(i\lambda_n\delta_1\cdots\delta_{\ell}\Theta)\ .
\]
The latter is then bounded by
\[
|{\rm Err}|\le x_n\sum_{\delta_1,\ldots,\delta_{\ell}=1}^{\infty}
\delta_{1}^{\ell-2}\delta_{2}^{\ell-3}\cdots\delta_{\ell-2}
\delta_{\ell-1}^{0}\delta_{\ell}^{-1}
\ e^{-\delta_1\cdots\delta_{\ell}t_n}\times
\frac{|i\lambda_n\delta_1\cdots\delta_{\ell}\Theta|^3}{6}\ ,
\]
i.e.,
\[
|{\rm Err}|\le
\frac{1}{6}x_n\Theta^3\lambda_n^3\ Z_{\ell+2}^{[\ell]}(t_n)\ .
\]
Using (\ref{xndefeq}), (\ref{lambdaasymeq}), (\ref{tnasymeq}), and
(\ref{Zplus2asymeq}), we see that the error is of order $n^{-\left(\frac{\ell-1}{2\ell}\right)}$ which goes to zero because of our standing assumption $\ell\ge 2$.

By our definition (\ref{lambdadefeq}) which is calibrated to keep the quadratic term invariant, the latter is exactly equal to $\frac{\Theta^2}{2}$.
On the other hand, the coefficient featuring in the linear term is small but not exactly zero. Namely, by our definitions (\ref{xndefeq}) and (\ref{tndefeq}), we have
\[
\lambda_n\left[n-x_n Z_{\ell}^{[\ell]}(t_n)\right]=\lambda_n \times n(1-e^{\frac{s}{b_n}})\ .
\]
By the hypothesis on the $b_n$ sequence in Proposition \ref{keyprop}, and (\ref{lambdaasymeq}), we get
\[
\lim\limits_{n\rightarrow\infty}
\lambda_n\left[n-x_n Z_{\ell}^{[\ell]}(t_n)\right]
=-s\sqrt{\ell-1}\ .
\]
We can finally use the dominated convergence theorem to deduce
\[
\lim\limits_{n\rightarrow\infty}
\int_{\mathbb{R}}f_n(\Theta)\ {\rm d}\Theta
=\int_{\mathbb{R}}e^{is\Theta\sqrt{\ell-1}-\frac{\Theta^2}{2}}\ {\rm d}\Theta
=\sqrt{2\pi}\ e^{-\frac{(\ell-1)s^2}{2}}\ ,
\]
by the familiar formula for the Fourier transform of the standard Gaussian.
The prefactor $\frac{\lambda_n}{2\pi}$ is a power law which still dominates the contribution of the minor arc, as seen from the bound (\ref{minacrusheq}), and therefore
\begin{equation}
\mathcal{J}_{\ell,n}(x_n,t_n)\sim
\frac{(x\mathcal{K}_{\ell})^{\frac{1}{2\ell}}}{\sqrt{2\pi\ell}}
\times n^{-\left(\frac{\ell+1}{2\ell}\right)}
\times e^{-\frac{(\ell-1)s^2}{2}}\ .
\label{Jasymeq}
\end{equation}

\section{Completion of the proofs of the main proposition and theorem}

We first expedite the proof of Proposition \ref{logHprop}.
By the decomposition (\ref{Hdecompeq}), the definition (\ref{Pdefeq}) of the prefactor $\mathcal{P}$, and the formula (\ref{Lniceeq}) for the function $\mathcal{L}$, we have
\[
\ln \mathcal{H}_{\ell,n}(x)=nt_n+xZ_{\ell-1}^{[\ell]}(t_n)+\ln \mathcal{J}_{\ell,n}(x,t_n)\ .
\] 
By (\ref{tnasymeq}), we obtain $nt_n=(x\mathcal{K}_{\ell})^{\frac{1}{\ell}}n^{\frac{\ell-1}{\ell}}(1+o(1))$, while using (\ref{tnasymeq}) as well as (\ref{Zminus1asymeq}), we also have
\begin{equation}
xZ_{\ell-1}^{[\ell]}(t_n)=
\frac{(x\mathcal{K}_{\ell})^{\frac{1}{\ell}}}{\ell-1}
n^{\frac{\ell-1}{\ell}}(1+o(1))\ .
\label{xZasymeq}
\end{equation}
The $s=0$ case of (\ref{Jasymeq}) shows that $\ln \mathcal{J}_{\ell,n}(x,t_n)$ is logarithmic and thus negligible compared to the two previous terms. Putting pieces together finishes proving Proposition \ref{logHprop}.

\medskip
Now we finish the proof of Proposition \ref{keyprop}.
By using (\ref{Jasymeq}) twice, for the possibly nonzero $s$, and for the $s=0$ special case, we see that the remainder term (\ref{remaindereq}) satisfies
\[
\lim\limits_{n\rightarrow\infty}R_{\ell,n}(s)=-\frac{(\ell-1)s^2}{2}\ .
\]
From (\ref{Psiconvenienteq}), we can write, by partially expanding the exponential,
\begin{eqnarray}
\Psi_{\ell,n}(s) &= & -\frac{s a_n}{b_n}+\left[nt_n+x_nZ_{\ell-1}^{[\ell]}(t_n) \right]
-\left[nt_n+xZ_{\ell-1}^{[\ell]}(t_n) \right]+R_{\ell,n}(s)\ \nonumber \\
 & = & -\frac{s a_n}{b_n}
+x(e^{\frac{s}{b_n}}-1)Z_{\ell-1}^{[\ell]}(t_n) +R_{\ell,n}(s) \nonumber \\
 & = & A_n s+ \frac{B_n s^2}{2}+x\left(e^{\frac{s}{b_n}}-1-\frac{s}{b_n}
-\frac{s^2}{2b_n^2}\right)Z_{\ell-1}^{[\ell]}(t_n) +R_{\ell,n}(s)\ ,
\label{lastPsieq}
\end{eqnarray}
where
\begin{eqnarray*}
A_n &:=& \frac{xZ_{\ell-1}^{[\ell]}(t_n)-a_n}{b_n}\\
B_n & := &\frac{xZ_{\ell-1}^{[\ell]}(t_n)}{b_n^2} \ .
\end{eqnarray*}
From the hypothese on the sequences from Proposition \ref{keyprop}, we immediately derive
\[
\lim\limits_{n\rightarrow\infty} A_n=0\ .
\]
From the same hypotheses, the part for $b_n$, and (\ref{xZasymeq}), we likewise derive
\[
\lim\limits_{n\rightarrow\infty} B_n=\ell\ ,
\]
as well as
\[
x\left(e^{\frac{s}{b_n}}-1-\frac{s}{b_n}
-\frac{s^2}{2b_n^2}\right)Z_{\ell-1}^{[\ell]}(t_n)=O\left(b_n^{-3}\right)\times
O\left(n^{\frac{\ell-1}{\ell}}\right)=O\left(n^{-\left(\frac{\ell-1}{2\ell}\right)}\right)
\]
which goes to zero. Mirroring the terms of (\ref{lastPsieq}), we see that
\[
\lim\limits_{n\rightarrow\infty}\Psi_{\ell,n}(s)=0+\frac{\ell s^2}{2}+0-\frac{(\ell-1)s^2}{2}=\frac{s^2}{2}\ ,
\]
and Proposition \ref{keyprop} is now established.

\medskip
To finish the proof of Theorem \ref{mainthm}, we first pick on purpose the sequences
\begin{eqnarray*}
a_n & :=& 
x Z_{\ell-1}^{[\ell]}(t_n) \\
b_n & := & 
\frac{(x\mathcal{K}_{\ell})^{\frac{1}{2\ell}}}{\sqrt{\ell(\ell-1)}}\times n^{\frac{\ell-1}{2\ell}}\ ,
\label{bnhypeq}
\end{eqnarray*}
which obviously satisfy the hypotheses of Proposition \ref{keyprop}. We then apply Theorem \ref{Curtissthm}, say with $s_0=1$, and collect the conclusion about the convergence of moments. In particular from the convergence of the first two moments, we readily obtain
\[
{\rm Var}(\mathsf{K}_{\ell,n})\sim b_n^2\sim \frac{(x\mathcal{K}_{\ell})^{\frac{1}{\ell}}}{\ell(\ell-1)}\times n^{\frac{\ell-1}{\ell}}
\]
and
\[
\mathbb{E}\mathsf{K_{\ell,n}}=x Z_{\ell-1}^{[\ell]}(t_n)+o\left(n^{\frac{\ell-1}{2\ell}}\right)
\]
which imply the asymptotics (\ref{Easymeq}) and (\ref{Vasymeq}).
Moreover, this also implies that if we now redefine
\begin{eqnarray*}
a_n & :=& \mathbb{E}\mathsf{K_{\ell,n}}
 \\
b_n & :=& \sqrt{{\rm Var}(\mathsf{K}_{\ell,n})} 
\ ,
\label{bnhypeq}
\end{eqnarray*}
these new sequences again satisfy the hypotheses of Proposition \ref{keyprop}. By another round of the latter followed by Theorem \ref{Curtissthm}, we conclude the proof of Theorem \ref{mainthm}, as stated.

\bigskip
\noindent{\bf Acknowledgements:}
{\small
A. A. thanks Bernhard Heim and Markus Neuhauser for instructive discussions related to the log-concavity of $\mathcal{H}_{\ell,n}(1)$ and related Bessenrodt-Ono inequalities.
A. A. also thanks Ofir Gorodetsky and Mark Wildon for enlighting comments or posts on MathOverflow related to the asymptotics of $Z_{1,0}^{[2]}(t)$. We again thank Ofir Gorodetsky for explaining to us the connection of our work to the CLT by Maples, Nikeghbali and Zeindler and the mean asymptotics by Ercolani and Ueltschi.
S. S. gratefully acknowledges a Simons collaboration grant.
S. S. also thanks Jim Pitman for pointing out useful references.
We are grateful to an anonymous referee for many useful suggestions including a suggestion for how to simplify the proof of Lemmas
2.2 and 2.3 and directing us to the reference \cite{PolyaSzego}.
} 

\appendix

\section{Generation of data for plots}
\label{sec:data}

In Figure \ref{fig:CLT} we compared the probability mass function for $\mathsf{K}_{\ell,n}$ to a Gaussian density function
with the same mean and variance.
We generated the data in two steps, which we fully explicate below.

\subsection{Wolfram Mathematica}
\small
\begin{verbatim}
DivisorSumList = Table[DivisorSum[n, # &],{n,1,150}];
BGenFun[z_] = Sum[DivisorSumList[[n]]/n*z^n,{n,1,150}];
DoubleGenFun[z_,x_] = Series[Exp[x*BGenFun[z]],{z,0,150}];
NumericalList = Log[CoefficientList[Expand[N[Coefficient[DoubleGenFun[z,x],z,150]]],x]];
\end{verbatim}
\normalsize
Then the output of NumericalList is a sequence of numbers that is the first command line for the Matlab file.

\subsection{Matlab file}
\small
\begin{verbatim}
Numerical_list = [0.908259,5.35459,8.96329,11.9068,14.3301,16.3321,17.9808,19.3252,...
    20.402,21.2396,21.8607,22.2839,22.5245,22.5957,22.5086,22.273,21.8972,...
    21.3887,20.7542,19.9996,19.1301,18.1507,17.0657,15.8791,14.5945,13.2152,...
    11.7445,10.1851,8.53962,6.81063,5.00036,3.11097,1.14444,-0.897335,-3.01261,...
    -5.19973,-7.45714,-9.78337,-12.177,-14.6368,-17.1615,-19.7499,-22.4009,...
    -25.1135,-27.8866,-30.7194,-33.6108,-36.5601,-39.5664,-42.6289,-45.747,...
    -48.9198,-52.1468,-55.4273,-58.7606,-62.1462,-65.5836,-69.0721,-72.6114,...
    -76.2009,-79.8401,-83.5287,-87.2661,-91.0521,-94.8862,-98.7681,-102.697,...
    -106.674,-110.697,-114.767,-118.883,-123.045,-127.253,-131.506,-135.805,...
    -140.149,-144.537,-148.971,-153.449,-157.972,-162.539,-167.15,-171.806,...
    -176.506,-181.25,-186.037,-190.869,-195.745,-200.664,-205.628,-210.635,...
    -215.687,-220.782,-225.921,-231.105,-236.333,-241.605,-246.922,-252.284,...
    -257.69,-263.141,-268.638,-274.18,-279.768,-285.401,-291.081,-296.808,...
    -302.581,-308.402,-314.27,-320.187,-326.152,-332.166,-338.23,-344.344,...
    -350.509,-356.725,-362.994,-369.316,-375.692,-382.122,-388.608,-395.151,...
    -401.752,-408.412,-415.133,-421.916,-428.762,-435.674,-442.653,-449.702,...
    -456.823,-464.019,-471.292,-478.646,-486.085,-493.613,-501.235,-508.957,...
    -516.786,-524.73,-532.799,-541.005,-549.364,-557.897,-566.632,-575.609,...
    -584.892,-594.6,-605.02]; 
Nmax=150;
kList = 1:Nmax;

figure 
hold on
axis([0,100,0.0001,0.17])

x_List = 0.1:0.1:Nmax;

color_List = ['krbg'];

for m=1:4
    r = 2*m-1;
    r_sqr = r^2;
    unnormalized_List = exp(Numerical_list + kList*log(r_sqr));
    prob_List = unnormalized_List/sum(unnormalized_List);
    mu_mean = sum(prob_List.*kList);
    sigma_var = sum(prob_List.*(kList-mu_mean).^2);
    log_z_List = -(x_List-mu_mean).^2/(2*sigma_var)-log(2*pi*sigma_var)/2;
    plot(kList,prob_List,[color_List(m),'.'],MarkerSize=15)
    plot(x_List,exp(log_z_List),"Color",color_List(m),LineWidth=2)
end\end{verbatim}


\begin{thebibliography}{999}

\bibitem{Abdesselam}
A. Abdesselam, Log-concavity with respect to the number of orbits for infinite tuples of commuting permutations. Ann. Comb. {\bf 29} (2025), no. 2, 563--573.
 
\bibitem{Abdesselam2025}
A. Abdesselam, Proof of a conjecture by Starr and log-concavity for random commuting permutations. Preprint arXiv:2506.06894[math.CO], 2025.

\bibitem{AbdesselamBDV}
A. Abdesselam, P. Brunialti, T. Doan, and P. Velie, 
A bijection for tuples of commuting permutations and a log-concavity conjecture.
Res. Number Theory {\bf 10} (2024), no. 2, Paper No. 45, 10 pp.

\bibitem{AhmadiGW}
L. Ahmadi, R. G\'{o}mez-A\'{\i}za,  and M. D. Ward,
A unified treatment of families of partition functions.
Matematica {\bf 3} (2024), no. 4, 1257--1296.

\bibitem{ArratiaBT}
R. Arratia, A. D. Barbour, and S. Tavar\'e,
Poisson process approximations for the Ewens sampling formula.
Ann. Appl. Probab. {\bf 2} (1992), no. 3, 519--535.

\bibitem{ArratiaBTbook}
R. Arratia, A. D. Barbour, and S. Tavar\'e,
{\it Logarithmic Combinatorial Structures: a Probabilistic Approach}.
EMS Monogr. Math.
European Mathematical Society (EMS), Z\"{u}rich, 2003.

\bibitem{BauerschmidtBS}
R. Bauerschmidt, D. C. Brydges, and G. Slade,
{\it Introduction to a Renormalisation Group Method}.
Lecture Notes in Math., {\bf 2242},
Springer, Singapore, 2019.

\bibitem{Billingsley}
P. Billingsley,
{\it Probability and Measure},
Second edition.
Wiley Ser. Probab. Math. Statist. Probab. Math. Statist.
John Wiley \& Sons, Inc., New York, 1986.

\bibitem{BridgesBBF}
W. Bridges, B. Brindle, K. Bringmann, and J. Franke,
Asymptotic expansions for partitions generated by infinite products.
Math. Ann. {\bf 390} (2024), no. 2, 2593--2632.

\bibitem{BringmannFH}
K. Bringmann, J. Franke, and B. Heim,
Asymptotics of commuting $\ell$-tuples in symmetric groups and log-concavity.
Res. Number Theory {\bf 10} (2024), no. 4, Paper No. 83, 19 pp.

\bibitem{deBruijn}
N. G. de Bruijn,
{\it Asymptotic Methods in Analysis},
Second edition.
Bibl. Math., Vol. {\bf IV},
North-Holland Publishing Co., Amsterdam; P. Noordhoff Ltd., Groningen, 1961. 

\bibitem{BryanF}
J. Bryan, and J. Fulman,
Orbifold Euler characteristics and the number of commuting $m$-tuples in the symmetric groups.
Ann. Comb. {\bf 2} (1998), no. 1, 1--6.

\bibitem{CantonFFM}
A. Cant\'on, J. Fern\'andez, P. Fern\'andez, and V. Maci\'a,
Khinchin families and Hayman class.
Comput. Methods Funct. Theory {\bf 21} (2021), no. 4, 851--904.

\bibitem{Crane}
H. Crane,
The ubiquitous Ewens sampling formula.
Statist. Sci. {\bf 31} (2016), no. 1, 1--19.

\bibitem{Curtiss}
J. H. Curtiss,
A note on the theory of moment generating functions.
Ann. Math. Statistics {\bf 13} (1942), 430--433.

\bibitem{ElboimG}
D. Elboim, and O. Gorodetsky,
Multiplicative arithmetic functions and the generalized Ewens measure.
Israel J. Math. {\bf 262} (2024), no. 1, 143--189.

\bibitem{ErcolaniU}
N. M. Ercolani, and D. Ueltschi,
Cycle structure of random permutations with cycle weights.
Random Structures Algorithms {\bf 44} (2014), no. 1, 109--133.

\bibitem{Ewens}
W. J. Ewens,
The sampling theory of selectively neutral alleles.
Theoret. Population Biol. {\bf 3} (1972), no. 1, 87--112.

\bibitem{Feller}
W. Feller,
The fundamental limit theorems in probability.
Bull. American Math. Soc. {\bf 51} (1945), 800--832.

\bibitem{FlajoletS}
P. Flajolet, and R. Sedgewick,
{\it Analytic combinatorics}.
Cambridge University Press, Cambridge, 2009. 

\bibitem{Goncharov1}
W. Gontcharoff,
Sur la distribution des cycles dans les permutations.
C. R. (Doklady) Acad. Sci. URSS (N.S.) {\bf 35} (1942), 267--269.

\bibitem{Goncharov2}
V. Gon\v{c}arov,
On the field of combinatory analysis.
American Math. Soc. Transl. (2) {\bf 19} (1962), 1--46.

\bibitem{GruberK}
C. Gruber, and H. Kunz,
General properties of polymer systems.
Comm. Math. Phys. {\bf 22} (1971), 133--161.

\bibitem{GuionnetMS}
A. Guionnet, and \'E. Maurel-Segala,
Matrix models at low temperature.
Preprint arXiv:2210.05239[math.PR], 2022. 

\bibitem{Hansen}
J. C. Hansen,
A functional central limit theorem for the Ewens sampling formula.
J. Appl. Probab. {\bf 27} (1990), no. 1, 28--43.

\bibitem{HeimN1}
B. Heim, and M. Neuhauser, Formulas for coefficients of polynomials assigned to arithmetic functions. Preprint arXiv:2010.07890[math.NT], (2020).

\bibitem{HeimN2}
B. Heim, and M. Neuhauser,
Horizontal and vertical log-concavity.
Res. Number Theory {\bf 7} (2021), no. 1, Paper No. 18, 12 pp.

\bibitem{KazakovKN}
V. A. Kazakov, I. K. Kostov, and N. Nekrasov,
$D$-particles, matrix integrals and KP hierarchy.
Nuclear Phys. B {\bf 557} (1999), no. 3, 413--442.

\bibitem{KazakovZ}
V. Kazakov, and Z. Zheng,
Analytic and numerical bootstrap for one-matrix model and ``unsolvable'' two-matrix model.
J. High Energy Phys. 2022, no. 6, Paper No. 30, 59 pp.

\bibitem{MageePvH}
M. Magee, D. Puder, and R. van Handel,
Strong convergence of uniformly random permutation representations of surface groups.
Preprint arXiv:2504.08988[math.GT], 2025.

\bibitem{MageeT} 
M. Magee, and J. Thomas,
Strongly convergent unitary representations of right-angled Artin groups.
Preprint arXiv:2308.00863[math.GR], 2023.

\bibitem{MaplesNZ}
K. Maples, A. Nikeghbali, and D. Zeindler,
On the number of cycles in a random permutation.
Electron. Commun. Probab. {\bf 17} (2012), no. 20, 13 pp.

\bibitem{McCarthy1}
J. E. McCarthy,
Random commuting matrices,
Preprint arXiv:2305.20029[math.PR], 2023.

\bibitem{McCarthy2}
J. E. McCarthy, and H. McCarthy,
Random anti-commuting Hermitian matrices.
Random Matrices Theory Appl. {\bf 13} (2024), no. 4, Paper No. 2450019, 13 pp.

\bibitem{Olshanski}
G. Olshanski,
Random permutations and related topics. In: {\it The Oxford Handbook of Random Matrix Theory}, edited by G. Akemann, J. Baik and P. Di Francesco,
pp. 510--533.
Oxford University Press, Oxford, 2011.

\bibitem{Pitman}
J. Pitman,
{\it Combinatorial stochastic processes}.
Lectures from the 32nd Summer School on Probability Theory held in Saint-Flour, July 7-24, 2002, with a foreword by J. Picard.
Lecture Notes in Math., {\bf 1875},
Springer-Verlag, Berlin, 2006.

\bibitem{PolyaSzego}
G. P\'olya and G. Szeg\"o,
{\it Problems and Theorems in Analysis, Volume I. Series. Integral Calculus. Theory of Functions. Translation by D. Aeppli.}
Springer-Verlag, Berlin, 1976.


\bibitem{Robin}
G. Robin,
Grandes valeurs de la fonction somme des diviseurs et hypoth\`ese de Riemann.
J. Math. Pures Appl. (9) {\bf 63} (1984), no. 2, 187--213.

\bibitem{StarrHR}
S. Starr, About the Hardy-Ramanujan partition function asymptotics.
Preprint arXiv:2408.08269[math.CO], 2024.

\bibitem{StarrAbund}
S. Starr, Some observations about the ``generalized abundancy index''.
Preprint arXiv:2505.07051[math.CO], 2025.

\bibitem{Tavare}
S. Tavar\'e,
The magical Ewens sampling formula.
Bull. London Math. Soc. {\bf 53} (2021), no. 6, 1563--1582.

\bibitem{Thomas}
A. Thomas,
Generalized punctual Hilbert schemes and $\mathfrak{g}$-complex structures.
Internat. J. Math. {\bf 33} (2022), no. 1, Paper No. 2250004, 44 pp.

\bibitem{Toth}
B. T\'oth,
Improved lower bound on the thermodynamic pressure of the spin 1/2 Heisenberg ferromagnet.
Lett. Math. Phys. {\bf 28} (1993), no. 1, 75--84.

\bibitem{Tripathi}
R. Tripathi,
On log-concavity of the number of orbits in commuting tuples of permutations.
Res. Number Theory {\bf 10} (2024), no. 4, Paper No. 78, 9 pp.

\bibitem{Zhang}
S. Zhang,
Log-concavity in powers of infinite series close to $(1-z)^{-1}$.
Res. Number Theory {\bf 8} (2022), no. 4, Paper No. 66, 17 pp.

\end{thebibliography}
\end{document}